\documentclass[11pt,oneside,english]{amsart}
\usepackage[T1]{fontenc}
\usepackage[latin9]{inputenc}
\usepackage{color}
\usepackage{multirow}
\usepackage{amsbsy}
\usepackage{amstext}
\usepackage{amsthm}
\usepackage{amssymb}
\usepackage{graphicx}
\usepackage{setspace}
\PassOptionsToPackage{normalem}{ulem}
\usepackage{ulem}
\makeatletter

\newcommand*\LyXbar{\rule[0.585ex]{1.2em}{0.25pt}}
\providecommand{\tabularnewline}{\\}
\numberwithin{equation}{section}
\numberwithin{figure}{section}
\theoremstyle{plain}
\newtheorem{thm}{\protect\theoremname}
\theoremstyle{plain}
\newtheorem{cor}[thm]{\protect\corollaryname}
\theoremstyle{remark}
\newtheorem{rem}[thm]{\protect\remarkname}

\usepackage{pdflscape}

\makeatother

\usepackage{babel}
\providecommand{\corollaryname}{Corollary}
\providecommand{\remarkname}{Remark}
\providecommand{\theoremname}{Theorem}

\begin{document}
\title{PDE Models and Riemann-Stieltjes Integrals in Sustainability}
\maketitle

\section*{Fulltitle: Partial Differential Equations Models and Riemann-Stieltjes
Integrals in Measuring Sustainability}
\begin{center}
\vspace{0.2cm}
\par\end{center}

\begin{center}
{\large{}Arni S.R. Srinivasa Rao{*}}{\large\par}
\par\end{center}

\begin{center}
Laboratory for Theory and Mathematical Modeling
\par\end{center}

\begin{center}
Medical College of Georgia
\par\end{center}

\begin{center}
Department of Mathematics, Augusta University,
\par\end{center}

\begin{center}
1120 15th Street, Augusta, GA 30912, USA
\par\end{center}

\begin{center}
Email: arrao@augusta.edu
\par\end{center}

\begin{center}
\vspace{0.05cm}
\par\end{center}

\begin{center}
\textbf{{*}} Corresponding author
\par\end{center}

\begin{center}
\vspace{0.3cm}
\par\end{center}

\begin{center}
{\large{}Sireesh Saride}{\large\par}
\par\end{center}

\begin{center}
Department of Civil Engineering,
\par\end{center}

\begin{center}
Indian Institute of Technology Hyderabad
\par\end{center}

\begin{center}
Kandi, Sangareddy - 502 285 Telangana, India
\par\end{center}

\begin{center}
Email: sireesh@ce.iith.ac.in
\par\end{center}

\vspace{0.5cm}

Appeared in \textbf{Advances on Methodology and Applications of Statistics (2021)
(Springer)} - 

A book in Honor of \textbf{C.R. Rao on the Occasion of his 100th Birthday
- edited by Barry C. Arnold, Narayanaswamy Balakrishnan, Carlos A.
Coelho. }

\vspace{0.2cm}

\begin{abstract}
Understanding sustainability through modeling involves one
of the complex and interdisciplinary activities where mathematics
plays a key role. We provide arguments favoring the need for developing
global models for measuring the status of sustainability. A global
model (applicable in broader perspective) and global sustainability
indices are proposed which can be used with real-world data. The solutions
of the proposed Partial Differential Equations (PDEs) are blended
with the weight functions of Riemann Stieltjes integrals to capture
the differential importance of sustainability associated factors.
The ideas, methods, and models are new and are prepared for handling
multi-dimensional and multi-variate data. A practically adaptable
formula for measuring the sustainability index is developed with few
key variables. We provide a real-world example arising in civil engineering
applications with a numerical example to demonstrate our models.
\end{abstract}


MSC: 92D40, 35Q80,26A42

\section{Introduction}

\label{sec:1} \textcolor{black}{According to the United States Environmental
Protection Agency \cite{USEPA} ``Sustainability creates and maintains
the conditions under which humans and nature can exist in productive
harmony, that permit fulfilling the social, economic and other requirements
of present and future generations'' and according to Oxford English
Dictionary \cite{OXD}, sustainable means - ``able to be maintained
at a certain rate or level''. Sustainability is a highly complex
and very broad-meaning word amalgamated from several cross-disciplinary
subjects and also influenced by political and governmental involvement.
Specific to Civil Engineering, American Society of Civil Engineers
(ASCE) defines sustainability as \textquotedbl a set of economic,
environmental and social conditions in which all of society has the
capacity and opportunity to maintain and improve its quality of life
indefinitely, without degrading the quantity, quality or the availability
of natural, economic and social resources (Ref: https://www.asce.org/advocacy/energy).
It is within and beyond the subject of academic collaboration. There
are discussions on achieving sustainability at academic and government
levels (see for example \cite{MAY}) and there is a need to address
sustainability from a common platform \cite{OST}. When we speak of
the mathematics of sustainability, the first question which comes
to me as a mathematical modeler is, can we build a global comprehensive
mathematical model to explain the status of sustainability? And secondly,
can we use it for predicting whether the system under investigation
is sustainable or not}\textbf{\textcolor{black}{?}}\textcolor{black}{{}
Supposing we can able to build one such model, then how does the mathematical
modeler (or mathematical community)}\textbf{\textcolor{black}{{}
}}\textcolor{black}{present it}\textbf{\textcolor{black}{{} }}\textcolor{black}{to
people with a non-mathematical background such that people who are
responsible for policy, government administrators and all other responsible
individuals understand the importance and implications of sustainability
and can foresee the predictions of a sustainability model? (It is,
of course, possible that people with a mathematical background occupy
government, political and other responsible professions). Mathematical
ideas cannot be adopted for the overall development of mankind and
for the betterment of species surrounded by mankind unless they can
be appreciated and encouraged by non-mathematicians. }

\textcolor{black}{A mathematician might develop an excellent model
to solve/explain certain practical issues, however, there are often
'gulfs' between `the people handling practical issues related to
sustainability' and `the mathematicians who are building and analyzing
models'. These gulfs maybe both ideological and conceptual. Ideological
gulfs can be broadly classified as undermining the applicability of
mathematical ideas and models for practical use in general, whereas
conceptual gulfs can be broadly classified as a lack of sufficient
background in mathematical sciences. Unless these concepts are accessible
to non-mathematicians, the fruits of mathematical reasoning would
not benefit mankind (and also surrounding species). But who will cross
the gap? Whether it is non-mathematicians who will cross the bridge
to reach mathematical thinkers or mathematicians who take their concepts
across the bridge could be a topic of consideration. One argument
could be that, since mathematicians produce mathematical models (methods),
they need to take them to other people and make them use models (like
a salesman sells his/her product). How often do mathematicians take
things in this direction? There are occasions when cross-disciplinary
teams including mathematicians build models for real-world solutions
in science and engineering. Although there may be an overlap in techniques,
modeling global sustainability is different from the modeling to explain
a particular phenomenon of science (or engineering), which is usually
conducted by a group of people with science (or engineering) and mathematical
backgrounds. There were several collaborative attempts with mathematicians
at building models for health policies that were conceived by people
responsible for policy formulation and successfully implemented at
the country level (for example, see \cite{ASR} and \cite{ASR2}). }

\textcolor{black}{Sustainability is a broader term than health policy
formulation. For example, if we have to formulate a mathematical model
based policy for controlling an epidemic and suppose we have the following
information: (i) health professionals who are handling population-level
treatment of this particular epidemic are aware of the transmission
dynamics of associated disease, (ii) time-series data on incidence
and prevalence, (iii) basic socio-demographic data for the population
of concern, (iv) intervention strategies by the government for controlling
the epidemic, (v) vaccine availability and distribution if at all
there exists a vaccination, (vi) any other relevant factors those
are essential in controlling the given epidemic. When a team of medical
and public health professionals with the above information approaches
a mathematician, the job of a mathematician is to understand the transmission
dynamics and develop a model for the same with all the variables and
parameters of interest. In case parameters are not readily available,
one has to use statistical procedures for estimation before model
building. It is also necessary to involve health professionals at
every stage of the model building for maximum flexibility and accuracy
of the model. Once a model is built and successfully tested for its
performance by taking inputs from health professionals in the team,
one can adapt the model for predicting the epidemic spread with and
without interventions. We can also measure the impact of certain policies
through such modeling. Understanding the sustainability status of
the population is not about to health alone. As mentioned previously,
it involves obtaining accurate information in all aspects of population
well-being, ranging from health, climate, food, agriculture resources,
science and technology to economic situation. Each category of this
information across various fields is required to use independently
and dependently in the models for sustainability status. Modeling
sustainability needs a cross-disciplinary team with an understanding
of cross-disciplinary data. See Figure schematic structure of the
process model building through the teamwork of cross-disciplinary
scientists, government agencies, politicians, etc,}

\begin{figure}
\includegraphics[scale=0.4]{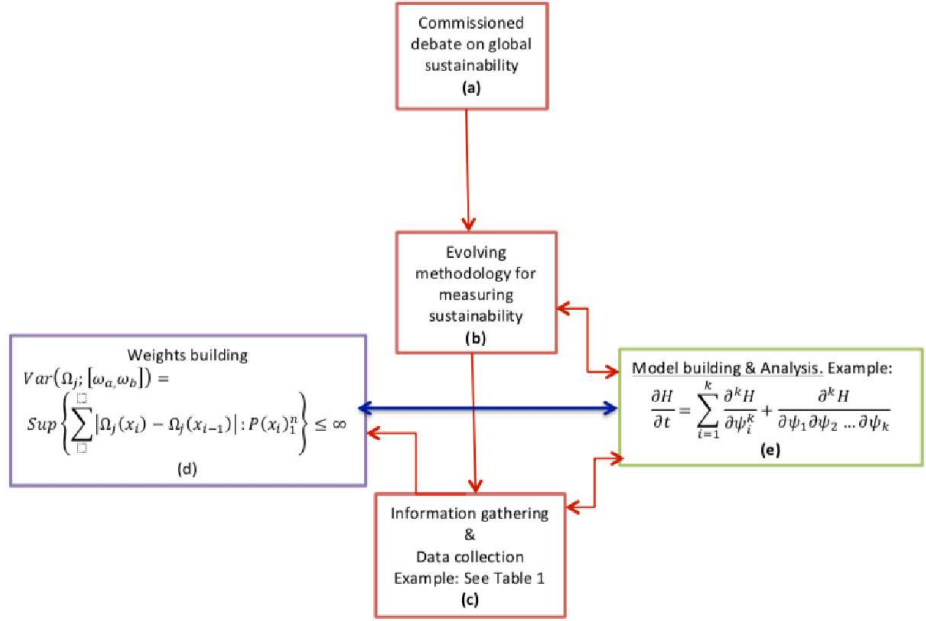}

\caption{\label{fig:Schematic-diagram-of}\textcolor{black}{Schematic diagram
of bridging between various sustainability stake holders. }\textbf{\textcolor{black}{(a)
}}\textcolor{black}{\uline{Commissioned debate on global sustainability}}\textcolor{black}{:
Team of scientists, politicians, government and non-government agencies
will debate on the status of sustainability, required further information,
direction of research needed for better understanding the status;
}\textbf{\textcolor{black}{(b) }}\textcolor{black}{\uline{Evolving
methodology for measuring sustainability}}\textcolor{black}{: Scientists
who are working on sustainability issues and organizations and individuals
who are working in practical implementation of maintaining sustainable
environment (see, for example, see www.epa.gov), statisticians and
mathematicians discuss as a team on the method of measuring sustainability
status; }\textbf{\textcolor{black}{(c) }}\textcolor{black}{\uline{Information
gathering \& data collection}}\textcolor{black}{: Required information
as decided in }\textbf{\textcolor{black}{(b) }}\textcolor{black}{will
be collected which will eventually help to arrive at Table 1; }\textbf{\textcolor{black}{(d)
}}\textcolor{black}{\uline{Weights building}}\textcolor{black}{:
With the help of }\textbf{\textcolor{black}{(c) }}\textcolor{black}{weight
functions will be formed and ordering of the weights by relative importance
in-terms of their contribution to the sustainability status is decided;
}\textbf{\textcolor{black}{(e) }}\textcolor{black}{\uline{Model
building \& analysis}}\textcolor{black}{: Models will be developed
since completion of }\textbf{\textcolor{black}{(b).}}\textcolor{black}{{}
The ideas generated in }\textbf{\textcolor{black}{(b), (c) }}\textcolor{black}{and
using the values obtained at }\textbf{\textcolor{black}{(d)}}\textcolor{black}{,
model based output will generated. There will be back-and-forth activities
between }\textbf{\textcolor{black}{(b), (c), (d) }}\textcolor{black}{and
}\textbf{\textcolor{black}{(e).}}}
\end{figure}

For measuring the status of sustainability, we need global models
with all-round global data, however, such a global model could have
components of sub-models (or local models) quantifying sustainability
status at various geographic regions on our planet. These computationally
intense models should be able to update global and local sustainability
status periodically such that time-dependent action-oriented policies
can be skimmed from these efforts. \textcolor{black}{Existing models
for sustainability involved variables from one or two categories and
these models are well constructed \cite{OST,Gund,Prett,Ravi}. Our
models proposed are new in terms of structure and conceptualization
and methods involve are different from the general sustainability
modeling framework proposed, see for example \cite{Philips,Todorov}.
We have proposed to use the weight functions of Riemann-Stieltjes
for giving differential importance to various factors involved in
sustainability measures.}\textbf{\textcolor{magenta}{{} }}Since the
length of the partitions in \textcolor{black}{Riemann-Stieltjes} could
be dynamically arranged, the corresponding weight functions proposed
in the work are flexible to capture variation in the sustainability
data of road transportation. In the next section, we propose sustainability
indices and models and show that these indices are solutions to proposed
models.

\section{Global Sustainability Models}

\label{sec:2} A general global model using \textcolor{black}{a Partial
Differential Equation (PDE)} for studying a measure of sustainability
(say, sustainability index, $H$ at time $t$), involving components
$\psi_{1},\psi_{2},...,\psi_{k}$ (i.e. independent variables) that
determine sustainability and partial derivatives of $H$ can be conceptualized
as one or more of the following equations:

\begin{eqnarray}
S\left(t,\psi_{1},\psi_{2},...,\psi_{k},H,H_{t},H_{\psi_{1}},H_{\psi_{2}},\cdots,H_{\psi_{k}}\right) & = & 0\label{eq:2.1.1}\\
\nonumber \\
S\left(\begin{array}{c}
t,\psi_{1},\psi_{2},...,\psi_{k},H,H_{t},H_{\psi_{1}},H_{\psi_{2}},\cdots,\\
H_{\psi_{k}},H_{\psi_{1}\psi_{2}},H_{\psi_{1}\psi_{3}},H_{\psi_{2}\psi_{3}},\cdots,H_{\psi_{1}\psi_{k}},H_{\psi_{2}\psi_{k}},\cdots
\end{array}\right) & = & 0\nonumber \\
\vdots\qquad\qquad &  & \vdots\label{eq:2.2.2}\\
S\left(\begin{array}{c}
t,\psi_{1},\psi_{2},...,\psi_{k},H,H_{t},H_{\psi_{1}},H_{\psi_{2}},\cdots,\\
H_{\psi_{k}},H_{\psi_{1}\psi_{2}},H_{\psi_{1}\psi_{3}},\cdots,H_{\psi_{2}\psi_{3,\cdots},\psi_{k}}
\end{array}\right) & = & 0\label{general-global-model}
\end{eqnarray}

where $H_{\psi_{1}},H_{\psi_{2}},...H_{\psi_{k}}$ are partial derivatives
and $H_{\psi_{1}\psi_{2}},H_{\psi_{1}\psi_{3}},H_{\psi_{2}\psi_{3}},$
$\cdots,H_{\psi_{1}\psi_{k}},H_{\psi_{2}\psi_{k}},\cdots$ $,H_{\psi_{2}\psi_{3,\cdots},\psi_{k}}$
are mixed partial derivatives.\textcolor{black}{{} For example}\textbf{\textcolor{black}{,}}

\textcolor{black}{{} 
\begin{eqnarray*}
H_{\psi_{i}} & = & \frac{\partial H}{\partial\psi_{i}},\\
H_{\psi_{1}\psi_{2}\cdots\psi_{k}} & = & \frac{\partial^{k}H}{\partial\psi_{1}\psi_{2}\cdots\psi_{k}},\\
H_{\psi_{k}\psi_{k}\cdots\psi_{k}} & = & \frac{\partial^{k}H}{\partial\psi_{k}^{k}},
\end{eqnarray*}
}

\textcolor{black}{and for}\textbf{\textcolor{black}{{} $\left|\alpha\right|=\alpha_{1}+\alpha_{2}+\cdots+\alpha_{k}$
}}\textcolor{black}{where $\alpha_{i}$ is a non-negative integer,
the mixed partial derivative, $H_{\psi_{k}\psi_{k}\cdots\psi_{k}}$,
means,}

\textcolor{black}{{} 
\begin{eqnarray*}
H_{\psi_{1}(\alpha_{1}\mbox{ times)}\psi_{k}(\alpha_{2}\mbox{ times)}\cdots\psi_{k}(\alpha_{k}\mbox{ times)}} & = & \frac{\partial^{k}H}{\partial\psi_{1}^{\alpha_{1}}\psi_{2}^{\alpha_{2}}\cdots\psi_{k}^{\alpha_{k}}}.
\end{eqnarray*}
}

\textcolor{black}{Here, $H:U\rightarrow\mathbb{R}$, where $U$ is
an open subset in $\mathbb{R}^{k}$ and $\psi_{i}\in U.$ (\ref{eq:2.1.1})
is first order, (\ref{eq:2.2.2}) is second order, and (\ref{general-global-model})
is $k^{th}$order PDE. In general,}

\textbf{\textcolor{black}{{} 
\begin{eqnarray}
S\left(t,\psi_{1},\psi_{2},...,\psi_{k},H,DH,D^{2}H,\cdots,D^{k-1}H,D^{k}H\right) & = & 0\mbox{ }\:\label{eq:general PDE}
\end{eqnarray}
}}

\textcolor{black}{for $\left(t,\psi_{1},\psi_{2},...,\psi_{k}\in U\right)$,
is called a $k^{th}$order PDE, where $S:\mathbb{R}^{k}\times\mathbb{R}^{k-1}\times\cdots\times\mathbb{R}\rightarrow\mathbb{R}$
and $D^{i}H$ is the set of all partial derivatives of order $i$
for $i=1,2,\cdots,k.$ The PDE (\ref{eq:general PDE}) could be linear
or non-linear. It is linear, by using $D^{\alpha}H$, and using two
functions, $f_{\alpha}\left(t,\psi_{1},\psi_{2},...,\psi_{k}\right)$
and $f\left(t,\psi_{1},\psi_{2},...,\psi_{k}\right)$, if we can write
in the form (\ref{eq:linear-form}), else, we can say it is not linear.}

\textbf{ 
\begin{eqnarray}
\sum_{\left|\alpha\right|\leq k}f_{\alpha}\left(t,\psi_{1},\psi_{2},...,\psi_{k}\right)D^{\alpha}H & = & f\left(t,\psi_{1},\psi_{2},...,\psi_{k}\right)\label{eq:linear-form}
\end{eqnarray}
}

The degree and combination of number of independent variables in $\Psi=\left(\psi_{1},\psi_{2},...,\psi_{k}\right)$
and construction of $H$ using $\Psi$ needs inputs from the global
community who are working on sustainability issues. The impact of
each of $\psi_{1},\psi_{2},...,\psi_{k}$ could be different in quantifying
$H$ and hence, we might have to introduce global weight functions
$\boldsymbol{\Omega}=\left(\Omega_{1},\Omega_{2},\cdots,\Omega_{k}\right)$
corresponding to each independent variable. We associate $\Omega$
with a Riemann-Stieltjes integrable function $F$ such that $F,\boldsymbol{\Omega}:[\omega_{a},\omega_{b}]\rightarrow\mathbb{R}$
are bounded functions on a compact interval $[\omega_{a},\omega_{b}]$.
A range of weights depending upon the fluctuations in independent
variables, for example, economic depression, shortage of resources
like food, energy, etc have to be assessed with cross-disciplinary
data such that a tagged partition function $P$ of $[\omega_{a},\omega_{b}]$
is obtained. Since $F$ is Riemann-Stieltjes integrable, we will have,

\begin{eqnarray}
\left|\Sigma\left(F,\Omega,P\right)-\int_{\omega_{a}}^{\omega_{b}}Fd\Omega\right| & < & \eta\label{RS}\\
\nonumber 
\end{eqnarray}

for every $\eta>0$. Here,

\begin{eqnarray}
\Sigma\left(F,\Omega,P\right) & = & \Sigma_{i=1}^{g}F(t_{i})\left[\Omega(x_{i})-\Omega(x_{i-1})\right]\label{RS-SUM}\\
\nonumber 
\end{eqnarray}

is the Riemann-Stieltjes sum of $F$ with respect to global sustainability
weight function $\Omega$ for $P=\left\{ \left(\left[x_{i-1},x_{i}\right],t_{i}\right)\right\} _{i=1}^{g}$,
\textcolor{black}{where $g$}\textbf{\textcolor{black}{{} }}\textcolor{black}{is
the size of partitions}\textbf{\textcolor{black}{. }}\textcolor{black}{For
definitions, properties and general description of Riemann-Stieltjes
integrals, refer to}\textbf{\textcolor{magenta}{{} \cite{Tao-book,Nil}}}\textbf{
}. Similarly, a general PDE for a specific geographic regional model
can be written in similar construction as that of model (\ref{general-global-model}).
The regional specific independent variables and regional specific
weights depending upon the regional sustainability indicators satisfying
the following Riemann-Stieltjes set up at the $h^{th}$ region are

\begin{eqnarray}
\left|\Sigma\left(F_{h},\Omega_{h},P_{h}\right)-\int_{\omega_{a}}^{\omega_{b}}F_{h}d\Omega_{h}\right| & < & \eta\label{eq:RS-REGIONAL}\\
\nonumber 
\end{eqnarray}

\begin{eqnarray}
\Sigma\left(F_{h},\Omega_{h},P_{h}\right) & = & \Sigma_{i=1}^{g_{1}}F_{h}(t_{i})\left[\Omega_{h}(x_{i})-\Omega_{h}(x_{i-1})\right]\label{eq:rsSUM-REGIONAL}\\
\nonumber 
\end{eqnarray}

\begin{eqnarray}
P_{h} & = & \left\{ \left(\left[x_{i-1},x_{i}\right],t_{i}\right)\right\} _{i=1}^{g_{1}}\label{partition-regional}\\
\nonumber 
\end{eqnarray}
\textcolor{black}{where $g_{1}$}\textbf{\textcolor{black}{{} }}\textcolor{black}{is
the size of partitions at the $h^{th}$ region.}\textbf{\textcolor{black}{{}
}}\textcolor{black}{Global models need to run simultaneously with
regional models over the time to obtain global sustainability measures.
Measurement of sustainability requires understanding influences of
more than one factor or variable and mutual inter-dependencies within
these variables. In fact, there will be several factors or variables
(see Table 1) that influence the population sustainability over time
period. Depending upon the situation one may consider two factors
or more factors together that influence the overall status of sustainability.
In Table 1, a list of variables provided for a general guidance. This
is a probable list and an effective list of sustainability influencing
could be emerged by splitting or combining variables from this list. }

\textcolor{black}{We propose two types of PDE models for $H$, one
without the mixed partial derivatives and second with the mixed partial
derivatives. First model (\ref{PDE-independent}) is a simple starting
point, because second order PDEs are known to have wider applications
in science and engineering (for example, see \cite{Evans-book,Strauss-book}).
This model is good when each dependent variable is considered in the
dynamics of sustainability status. When we need to incorporate interaction
(or influence) of more variables then model (\ref{PDE-dependent})
will be suitable. In this model we have considered the term }\textbf{\textcolor{black}{$\frac{\partial^{k}H}{\partial\psi_{1}\cdots\partial\psi_{k}}$}}\textcolor{black}{{}
for demonstration of our analysis. One can consider other mixed partial
derivatives such as$\mathnormal{\frac{\partial^{k}H}{\partial\psi_{3}\partial\psi_{2}\partial\psi_{1}\cdots\partial\psi_{k}}}\,$
or other combination or order of the partial derivatives. For demonstration
of the solution in this paper, we use the model (\ref{PDE-dependent}),
and similar approach can be adapted for all other models with different
mixed partial derivatives, for example see (\ref{PDE-dependent})
and (\ref{PDE-dependent}).}

\begin{eqnarray}
\frac{\partial H}{\partial t} & = & \Sigma_{i=1}^{k}\frac{\partial^{2}H}{\partial\psi_{i}^{2}}\,\left(k>1\right)\label{PDE-independent}\\
\frac{\partial H}{\partial t} & = & \Sigma_{i=1}^{k}\frac{\partial^{k}H}{\partial\psi_{i}^{k}}+\frac{\partial^{k}H}{\partial\psi_{1}\cdots\partial\psi_{k}}\,\left(k>1\right)\label{PDE-dependent}
\end{eqnarray}

If we ignore mixed partial derivatives in Model (\ref{PDE-dependent})
and fix $k=2$, then model (\ref{PDE-dependent}) becomes model (\ref{PDE-independent}).
\textcolor{black}{Other PDE model by varying mixed partial derivatives
as described above are, }

\textcolor{black}{{} 
\begin{eqnarray}
\frac{\partial H}{\partial t} & = & \Sigma_{i=1}^{k}\frac{\partial^{k}H}{\partial\psi_{i}^{k}}+\frac{\partial^{k}H}{\partial\psi_{3}\partial\psi_{2}\partial\psi_{1}\cdots\partial\psi_{k}}\,\left(k>1\right)\nonumber \\
{\color{black}} & {\color{black}} & {\color{black}}\\
\frac{\partial H}{\partial t} & = & {\color{black}\Sigma_{i=1}^{k}\frac{\partial^{k}H}{\partial\psi_{i}^{k}}+\frac{\partial^{k}H}{\partial\psi_{1}\cdots\partial\psi_{k}}+\frac{\partial^{k-1}H}{\partial\psi_{1}\cdots\partial\psi_{k-1}}+\,\,\left(k>1\right)}\nonumber \\
{\color{black}} & {\color{black}} & {\color{black}\frac{\partial^{k-2}H}{\partial\psi_{3}\partial\psi_{2}\partial\psi_{1}\cdots\partial\psi_{k-2}}+\cdots+\frac{\partial H}{\partial\psi_{k}}\,}
\end{eqnarray}
}

\textcolor{black}{Actual function form of a PDE model with a suitable
order and degree are usually decided by the experts working in the
field and in this case the team of cross-disciplinary scientists after
sufficient debate. We confine to the model (\ref{PDE-dependent})
for obtaining a sustainability index to be used by groups or individuals
working on practical modeling. } 
\begin{thm}
\label{thm1}For the Model (\ref{PDE-independent}), $(i)H=kt+\frac{1}{2}\Sigma_{i=1}^{k}\psi_{i}^{2}$
and $(ii)H=(2k)t+\Sigma_{i=1}^{k}\psi_{i}^{2}$ are solutions. 
\end{thm}

\begin{proof}
(i) For $k=1$, the model (\ref{PDE-independent}) is a heat equation
and the solution is $H=t+\frac{1}{2}\psi_{1}^{2}$. For $H=kt+\frac{1}{2}\Sigma_{i=1}^{k}\psi_{i}^{2}$,
we have $\frac{\partial H}{\partial t}=k=\Sigma_{i=1}^{k}\frac{\partial^{2}H}{\partial\psi_{i}^{2}}.$

(ii) For $H=(2k)t+\Sigma_{i=1}^{k}\psi_{i}^{2}$ we have $\frac{\partial H}{\partial t}=2k=\Sigma_{i=1}^{k}\frac{\partial^{2}H}{\partial\psi_{i}^{2}}$.
Hence the theorem is proved. 
\end{proof}
\begin{thm}
\label{thm2}For the Model (\ref{PDE-dependent}), (i) $H=(k+1)t+\frac{1}{k!}\Sigma_{i=1}^{k}\psi_{i}^{k}+\Pi_{i=1}^{k}\psi_{i}$
and

(ii) $H=\left(\left(k!k\right)+1\right)t+\Sigma_{i=1}^{k}\psi_{i}^{k}+\Pi_{i=1}^{k}\psi_{i}$
are solutions. 
\end{thm}

\begin{proof}
(i) Let $H=(k+1)t+\frac{1}{k!}\Sigma_{i=1}^{k}\psi_{i}^{k}+\Pi_{i=1}^{k}\psi_{i}$.
We have $\frac{\partial H}{\partial\psi_{i}}=\frac{k}{k!}\psi_{i}^{k-1}+\psi_{1}\cdots\psi_{i-1}\psi_{i+1}\cdots\psi_{k}$
and $\frac{\partial^{k}H}{\partial^{k}\psi_{i}}=1$. We have, 
\begin{eqnarray*}
\frac{\partial^{k}H}{\partial\psi_{1}\cdots\partial\psi_{k}} & = & \frac{\partial^{k-1}}{\partial\psi_{1}\cdots\partial\psi_{k-1}}\left(\frac{\partial H}{\partial\psi_{k}}\right)\\
 & = & \frac{\partial^{k-1}}{\partial\psi_{1}\cdots\partial\psi_{k-1}}\left(\frac{k}{k!}\psi_{k}^{k-1}+\Pi_{i=1}^{k-1}\psi_{i}\right)\\
 & = & \frac{\partial^{k-2}}{\partial\psi_{1}\cdots\partial\psi_{k-2}}\left(\Pi_{i=1}^{k-2}\psi_{i}\right)=1
\end{eqnarray*}

Therefore,

$\frac{\partial H}{\partial t}=k+1=\Sigma_{i=1}^{k}\frac{\partial^{k}H}{\partial\psi_{i}^{k}}+\frac{\partial^{k}H}{\partial\psi_{1}\cdots\partial\psi_{k}}\,\left(k>1\right)$,
which confirms that given $H$ is a solution.

(ii) Let $H=\left(\left(k!k\right)+1\right)t+\Sigma_{i=1}^{k}\psi_{i}^{k}+\Pi_{i=1}^{k}\psi_{i}$.
We have $\frac{\partial H}{\partial\psi_{i}}=k\psi_{i}^{k-1}+\psi_{1}\cdots\psi_{i-1}\psi_{i+1}\cdots\psi_{k}$
and $\frac{\partial^{k}H}{\partial^{k}\psi_{i}}=k!$. We have,

\begin{eqnarray*}
\frac{\partial^{k}H}{\partial\psi_{1}\cdots\partial\psi_{k}} & = & \frac{\partial^{k-1}}{\partial\psi_{1}\cdots\partial\psi_{k-1}}\left(\frac{\partial H}{\partial\psi_{k}}\right)\\
 & = & \frac{\partial^{k-1}}{\partial\psi_{1}\cdots\partial\psi_{k-1}}\left(k\psi_{k}^{k-1}+\Pi_{i=1}^{k-1}\psi_{i}\right)\\
 & = & \frac{\partial^{k-2}}{\partial\psi_{1}\cdots\partial\psi_{k-2}}\left(\Pi_{i=1}^{k-2}\psi_{i}\right)=1
\end{eqnarray*}

Therefore,

$\frac{\partial H}{\partial t}=k!k+1=\Sigma_{i=1}^{k}\frac{\partial^{k}H}{\partial\psi_{i}^{k}}+\frac{\partial^{k}H}{\partial\psi_{1}\cdots\partial\psi_{k}}\,\left(k>1\right)$,
which confirms that given $H$ is a solution. 
\end{proof}
\begin{cor}
Suppose $\alpha>0$ and $\beta>0$ are two parameters, then the function,
$H(t,\psi_{1},\cdots,\psi_{k};\alpha,\beta)=\left(k\alpha k!+\beta\right)t+\frac{1}{k!}\Sigma_{i=1}^{k}\psi_{i}^{k}+\Pi_{i=1}^{k}\psi_{i}$
is also a solution of the model (\ref{PDE-dependent}). 
\end{cor}

These two theorems provide basic general idea of solutions for the
models in (\ref{PDE-independent}) and (\ref{PDE-dependent}). These
solutions are two proposed candidates for sustainability indices without
weights. In the next theorem we provide a solution for above models
(we will call weighted sustainability index, $H(\Omega$)), which
has weight functions described in this section. 
\begin{thm}
\label{thm:The-weighted-sustainability index}The weighted sustainability
index

$H(\Omega)=\left(k!\Sigma_{i=1}^{k}\Omega_{i}+\prod_{i=1}^{k}\Omega_{i}\right)t$
$+\Sigma_{i=1}^{k}\Omega_{i}\psi_{i}^{k}+\prod_{i=1}^{k}\Omega_{i}\psi_{i}$
is a solution for the sustainability model (\ref{PDE-dependent})
with mixed partial derivatives. 
\end{thm}

\begin{proof}
We have, $\frac{\partial H(\Omega)}{\partial\psi_{i}}=\Omega_{i}k\psi_{i}^{k-1}+\Omega_{1}\psi_{1}\cdots\Omega_{i-1}\psi_{i-1}\Omega_{i}\Omega_{i+1}\psi_{i+1}\cdots\Omega_{k}\psi_{k}$
and $\frac{\partial^{k}H(\Omega)}{\partial^{k}\psi_{i}}=\Omega_{i}k!$.
The mixed partial derivative terms in the model (\ref{PDE-dependent})
can be obtained for the given index as,

\begin{eqnarray*}
\frac{\partial^{k}H(\Omega)}{\partial\psi_{1}\cdots\partial\psi_{k}} & = & \frac{\partial^{k-1}}{\partial\psi_{1}\cdots\partial\psi_{k-1}}\left(\frac{\partial H(\Omega)}{\partial\psi_{k}}\right)\\
 & = & \frac{\partial^{k-1}}{\partial\psi_{1}\cdots\partial\psi_{k-1}}\left(\Omega_{k}k\psi_{k}^{k-1}+\Omega_{k}\Pi_{i=1}^{k-1}\psi_{i}\right)\\
 & = & \frac{\partial^{k-2}}{\partial\psi_{1}\cdots\partial\psi_{k-2}}\left(\Omega_{k}\Omega_{k-1}\Pi_{i=1}^{k-2}\psi_{i}\right)\\
 & = & \prod_{i=1}^{k}\Omega_{i}
\end{eqnarray*}

Therefore,

\begin{eqnarray*}
\frac{\partial H(\Omega)}{\partial t} & =k!\Sigma_{i=1}^{k}\Omega_{i}\psi_{i}+\prod_{i=1}^{k}\Omega_{i}= & \Sigma_{i=1}^{k}\frac{\partial^{k}H}{\partial\psi_{i}^{k}}+\frac{\partial^{k}H}{\partial\psi_{1}\cdots\partial\psi_{k}}\,\left(k>1\right)\\
\end{eqnarray*}

Hence, the weighted function proposed is a solution of the model (\ref{PDE-dependent}). 
\end{proof}
\begin{cor}
\label{cor:1} From theorem (\ref{thm:The-weighted-sustainability index}),
we can see that

\textup{$H(\Omega)=\left(\Sigma_{i=1}^{k}\Omega_{i}+\prod_{i=1}^{k}\Omega_{i}\right)t+\frac{1}{k!}\Sigma_{i=1}^{k}\Omega_{i}\psi_{i}^{k}+\prod_{i=1}^{k}\Omega_{i}\psi_{i}$
is also a solution for the model (\ref{PDE-dependent}). } 
\end{cor}

$ $ 
\begin{cor}
\label{cor2, alpha, beta}Suppose $\alpha>0$ and $\beta>0$ are two
parameters, then the function, 
\begin{eqnarray*}
H(t,\psi_{1},\cdots,\psi_{k},\Omega_{1},\cdots,\Omega_{k};\alpha,\beta) & = & \left(\alpha k!\Sigma_{i=1}^{k}\Omega_{i}+\beta\prod_{i=1}^{k}\Omega_{i}\right)t\\
 &  & \quad+\alpha\Sigma_{i=1}^{k}\Omega_{i}\psi_{i}^{k}+\beta\prod_{i=1}^{k}\Omega_{i}\psi_{i}
\end{eqnarray*}
is also a solution of the model (\ref{PDE-dependent}). 
\end{cor}

\textcolor{black}{We shortlist seven key independent variables based
on \cite{MAY} and also using self-intuition in measuring sustainability.}\textbf{\textcolor{black}{{}
}}\textcolor{black}{Seven variables are, food and agriculture $(\psi_{1})$,
climate and environment $(\psi_{2})$, population and economics $(\psi_{3})$,
political situation}\textbf{\textcolor{black}{{} $(\psi_{4})$}}\textcolor{black}{,
medical technology $(\psi_{5})$, energy $(\psi_{6})$, and science
and technology $(\psi_{7})$. These variables reduce the models (\ref{PDE-independent})
and (\ref{PDE-dependent}) with a lesser number of realistic variables.
Models with reduced variables are to be run along with regional level
sub-models with regional specific weight functions and regional specific
combinations of key variables listed. A hypothetical description of
these seven variables and expected data at a country or a concerned
region is given in Table \ref{tab:Seven-variables-required}.}\textbf{\textcolor{black}{{}
}}\textcolor{black}{Some of these variables can be split into two
or more variables if there is enough evidence from the data. For smaller
regions, one can have a larger number of variables because the data
needed for obtaining an index would be relatively easier to collect.
For larger countries, the more the number of variables more will be
the potential variation in the index.}\textbf{\textcolor{black}{{}
}}\textcolor{black}{We suggest $H(\Omega_{7})$ and $H(\Omega_{7};\alpha,\beta)$
as two candidates for the sustainability index using the }\textbf{\textcolor{black}{seven}}\textcolor{black}{{}
variables. }

The\textcolor{black}{{} seven variable sustainability index based
on corollary (\ref{cor:1}) is,}

\textcolor{black}{{} 
\begin{eqnarray}
H(\Omega_{7}) & = & \left(\Sigma_{i=1}^{7}\Omega_{i}+\prod_{i=1}^{7}\Omega_{i}\right)t+\frac{1}{7!}\Sigma_{i=1}^{7}\Omega_{i}\psi_{i}^{7}+\prod_{i=1}^{7}\Omega_{i}\psi_{i}\label{Index with 6 variables}\\
\nonumber 
\end{eqnarray}
}

\textcolor{black}{and the seven variable sustainability index based
on corollary (\ref{cor2, alpha, beta}) is,}

\textcolor{black}{{} 
\begin{eqnarray}
H(\Omega_{7};\alpha,\beta) & = & \left(\alpha6!\Sigma_{i=1}^{7}\Omega_{i}+\beta\prod_{i=1}^{7}\Omega_{i}\right)t\nonumber \\
 &  & \quad+\alpha\Sigma_{i=1}^{6}\Omega_{i}\psi_{i}^{7}+\beta\prod_{i=1}^{7}\Omega_{i}\psi_{i}\label{eindex with six var, alpha, beta}
\end{eqnarray}
}

\textcolor{black}{Once we have data for the variables in Table \ref{tab:Seven-variables-required}
are collected and weights are computed, the index given in (\ref{eindex with six var, alpha, beta})
or similar index with different mixed partial derivatives can be constructed.
We have two more parameters, $\alpha$ and $\beta$ in (\ref{eindex with six var, alpha, beta})
which can be fitted using a least square method or suitable method
depending upon the dimension of the data. Weight computation is needed
to be done with the help of experts in the respective field. The relative
importance of various weights can be decided by the entire team to
capture the differential influences of variables in global sustainability
status.}

\begin{landscape}

\begin{table}
\begin{tabular}{|c|c|c|}
\hline 
\textcolor{black}{Variable} & \textcolor{black}{Description} & \textcolor{black}{Symbol}\tabularnewline
\hline 
\hline 
\textcolor{black}{food and agriculture} & \textcolor{black}{proportion of yield of food production to the arable
land} & \textcolor{black}{$\psi_{1}$}\tabularnewline
\hline 
\textcolor{black}{climate and environment} & \textcolor{black}{$\begin{array}{c}
\text{tproportion of people who are living under the clean air,}\\
\text{ clearn water facilities and normal climatic conditions}
\end{array}$} & \textcolor{black}{$\psi_{2}$}\tabularnewline
\hline 
\textcolor{black}{population and economics} & \textcolor{black}{$\begin{array}{c}
\text{proportion of people who are above he poverty levels }\\
\text{prescribed by the multidimensional poverty index}\\
\begin{array}{c}
\text{and proportion of skilled population required for }\\
\text{economic growth}
\end{array}
\end{array}$} & \textcolor{black}{$\psi_{3}$}\tabularnewline
\hline 
\textcolor{black}{political situation} & \textcolor{black}{$\begin{array}{c}
\text{duration of the time period where the population }\\
\text{concerned or a region or a country is living under }\\
\text{stable political situation as per global perspective }
\end{array}$} & \textcolor{black}{$\psi_{4}$}\tabularnewline
\hline 
\textcolor{black}{medical technology} & $\begin{array}{c}
\text{proportion of people who are living}\\
\begin{array}{c}
\text{within the reach of best medical treatment }\\
\text{and health facilities}
\end{array}
\end{array}$ & \textcolor{black}{$\psi_{5}$}\tabularnewline
\hline 
\textcolor{black}{energy} & \textcolor{black}{$\begin{array}{c}
\text{proportion of industrial sector, transportation sector}\\
\text{and agricultural sector, general population infrastructure }\\
\begin{array}{c}
\text{which are availaling sufficient electricity for }\\
\text{optimum productivity }
\end{array}
\end{array}$} & \textcolor{black}{$\psi_{6}$}\tabularnewline
\hline 
\textcolor{black}{science and technology} & \textcolor{black}{$\begin{array}{c}
\text{level of research in basic sciences and}\\
\text{overall technological advancements or }\\
\text{availability of technology to support the needs }\\
\text{of concerned population or a region or a country}
\end{array}$} & \textcolor{black}{$\psi_{7}$}\tabularnewline
\hline 
\end{tabular}

\vspace{0.2cm}

\caption{\label{tab:Seven-variables-required}\textcolor{black}{Seven variables
required for the sustainability index}}
\end{table}

\end{landscape}

\section{Properties of $H\left(\Omega\right)$}

\label{sec:3} We investigate properties of the sustainability index
proposed in section \ref{sec:2}. For measuring $H(\Omega_{j})$ we
need practical weight functions which is possible by assuming

\begin{eqnarray}
\Sigma_{i=1}^{n}\left|\Omega_{j}(x_{i})-\Omega_{j}(x_{i-1})\right| & \leq & M\in\mathbb{R}\mbox{ for each }j=1,2,\cdots,k\label{eq:BV(OMEGA)}\\
\nonumber 
\end{eqnarray}

for all partitions of $\left[\omega_{a},\omega_{b}\right].$\textcolor{black}{{}
}\textbf{\textcolor{black}{I}}\textcolor{black}{n this case, $\Omega_{j}$
is an increasing function. Since $\Omega_{j}:[\omega_{a},\omega_{b}]\rightarrow\mathbb{R},$
the variation of $\Omega_{j}$ over $[\omega_{a},\omega_{b}]$ is }

\textcolor{black}{{} 
\begin{eqnarray*}
Var\left(\Omega_{j};[\omega_{a},\omega_{b}]\right) & = & \sup\left\{ \Sigma_{i=1}^{n}\left|\Omega_{j}(x_{i})-\Omega_{j}(x_{i-1})\right|:P\right\} \leq\infty
\end{eqnarray*}
This implies, $\Omega_{j}$ has bounded variation (BV) on $\left[\omega_{a},\omega_{b}\right]$
and written as, }\textbf{\textcolor{magenta}{$ $}}$\Omega_{j}\in BV\left(\left[\omega_{a},\omega_{b}\right]\right)$
and $\Omega_{j}$ is bounded on $\left[\omega_{a},\omega_{b}\right]$.
Each of the $\Omega_{j}$ is constructed from the $\Omega_{j}(h)$,
which is a weight function for $h^{th}-$region. 
\begin{thm}
For an $\Omega_{j}$ satisfying (\ref{eq:BV(OMEGA)}), 
\begin{eqnarray*}
\sup\left\{ \Sigma_{i=1}^{n}\left|\Omega_{j}(x_{i})-\Omega_{j}(x_{i-1})\right|:P\right\}  & \geq & \frac{\left|\int_{\omega_{a}}^{\omega_{b}}Fd\Omega_{j}\right|}{\left[\sup_{t\in\left[\omega_{a},\omega_{b}\right]}\left|F(t)\right|\right]}.\\
\end{eqnarray*}
\end{thm}

\begin{proof}
From the existence theorem for the Riemann-Stieltjes integral (see
for example, \cite{Nil}), we have,

\begin{eqnarray*}
\left|\int_{\omega_{a}}^{\omega_{b}}Fd\Omega_{j}\right| & \leq & \left[\sup_{t\in\left[\omega_{a},\omega_{b}\right]}\left|F(t)\right|\right]\times\\
 &  & \quad Var\left(\Omega_{j};\left[\omega_{a},\omega_{b}\right]\right)\\
\implies\left|\int_{\omega_{a}}^{\omega_{b}}Fd\Omega_{j}\right| & \leq & \left[\sup_{t\in\left[\omega_{a},\omega_{b}\right]}\left|F(t)\right|\right]\times\\
 &  & \sup\left\{ \Sigma_{i=1}^{n}\left|\Omega_{j}(x_{i})-\Omega_{j}(x_{i-1})\right|:P\right\} \\
 &  & \left(\mbox{ for all partitions }P\mbox{ of }\left[\omega_{a},\omega_{b}\right]\right)\\
\implies\sup\left\{ \Sigma_{i=1}^{n}\left|\Omega_{j}(x_{i})-\Omega_{j}(x_{i-1})\right|:P\right\}  & \geq & \frac{\left|\int_{\omega_{a}}^{\omega_{b}}Fd\Omega_{j}\right|}{\left[\sup_{t\in\left[\omega_{a},\omega_{b}\right]}\left|F(t)\right|\right]}\\
\end{eqnarray*}
\end{proof}
\begin{rem}
Suppose $\Omega_{j}(\omega_{a})=A_{j1}$ and $\Omega_{j}(\omega_{b})=A_{j2}$,
and assuming $\Omega_{j}$ is an $1-1$ function for all $j=1,2,\cdots,k$,
\textcolor{black}{then using Theorem}\textbf{\textcolor{black}{{}
\ref{thm:The-weighted-sustainability index} }}\textcolor{black}{and
Corollary}\textbf{\textcolor{black}{{} \ref{cor:1}, }}\textcolor{black}{we
have, $\frac{\partial H(\Omega_{j})}{\partial t}=\left(\Sigma_{j=1}^{k}\left[A_{j1},A_{j2}\right]+\prod_{j=1}^{k}\left[A_{j1},A_{j2}\right]\right)$.} 
\end{rem}

\section{Real World Situations: Recycling in Pavement Construction}

\label{sec:4}
\begin{figure}
\begin{enumerate}
\item \includegraphics[scale=0.5]{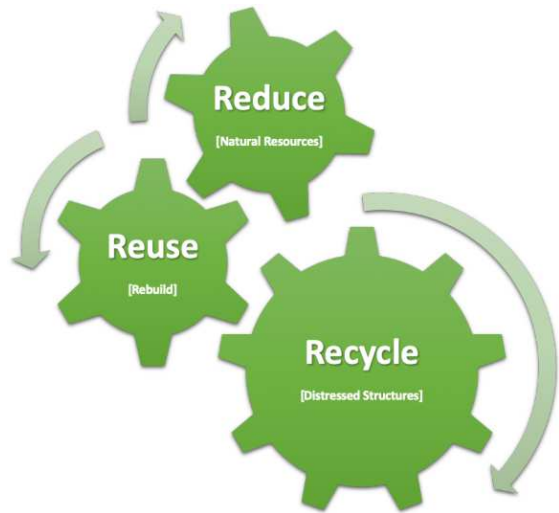}
\end{enumerate}
\caption{\label{fig:Sustainable-life-cycle}Sustainable life cycle of a construction
materials -- When the civil structures are demolished and recycles
at the end of their design life, can be recycled and prepare for reuse
them as a resource materials in a new and rehabilitation works. When
we reuse the reclaimed materials, the demand on the natural materials
can be reduced. This cycle continues to reduce the use of natural
resources and increases the sustainability index.}
\end{figure}

The availability of high-quality aggregates for the
construction of highways is scarce across the world; hence the cost
of construction has become prohibitively high due to increased lead
distance to mobilize suitable quality materials. A simple calculation
shows that approx. Twelve thousand tons of Natural Aggregates are
required per lane kilometer to build an express highway. If locally
available marginal or recycled aggregates can be used after a proper
stabilization process can solve the issue. Besides, recycling the
existing distressed pavements and reusing them back in the construction
of new and rehabilitated pavements has become a viable alternative
to natural aggregates. Recycling can address some of the sustainability
issues such as the conservation of natural resources and fossil fuels,
preservation of the environment, and retention of the existing highway
geometrics \cite{Taha A}. Reusing the construction and demolition
waste generated from the construction industry can reduce the burden
on the landfills, thus reducing the impact on the environment; otherwise,
the waste would have been ended up in landfills. The sustainable cycle
of the construction industry can be seen in Figure \ref{fig:Sustainable-life-cycle}.

There are different types of recycled and/or reclaimed materials generated
across the world from different processes, include but not limited
to construction and demolition waste (C\&D), crushed bricks (CB),
recycled concrete aggregate (RCA), reclaimed asphalt pavement (RAP),
quarry waste, (QA), recycled glass (RG), roofing shingles (RS), etc.
Among these recycled materials, the RAP, RAC, and C\&D are generated
in substantial quantities. The Federal Highway Administration (FHWA)
has estimated that close to 100 million tons of RAP material is produced
by milling HMA (Hot mix asphalt) each year \cite{MO Ashpalt}. According
to United States Geological Survey (USGS), the use of recycled crushed
concrete aggregates is still constituted by only about 0.5\% of the
total aggregates consumed in the US \cite{Gnanendran}. As high as
41 million tons of quarry by-products are generated per annum in the
United Kingdom \cite{Manning}.

However, these recycled aggregates `as is' can't meet the structural
and strength characteristics required for various civil engineering
applications, including pavement construction due to their inferior
characteristics. Several researchers have pointed that stabilizing
RAP and natural virgin aggregate (VA) mixes with conventional cementitious
materials like cement and lime have yielded superior strength and
stiffness properties \cite{Taha A,Puppala-saride}. However, they
could not promote the high replacement of VA with RAP. In the US,
state transportation agencies allowed only up to 20\% RAP due to a
lack of understanding of their behavior. Nevertheless, the production
of Portland cement is associated with a substantial amount of energy
depletion (5,000 MJ/ton cement), non-renewable resources (1.5-ton
limestone and clay/ton cement) as well as $CO_{2}$ emissions (0.95
ton $CO_{2}$/ton cement) which leads to climate change and ecological
imbalance \cite{Higgins,YiYLisks}. Hence, there is a necessity to
ascertain more environmentally friendly and green materials that can
match at least the performance of the traditional materials (cement
and lime) under similar circumstances. Thus the utilization of potential
industrial by-products is envisioned to replace either in the partial
or full amount of conventional stabilizers.

With the substantial availability of industrial by-products such as
fly ash and ground granulated blast furnace slag (GGBS), stabilizing
RAP:VA blends with these materials could reduce not only the cost
of construction but also protect the environment \cite{SarideAvimeni,ArulRajah,SarideJallu}.
Besides, Agarwal et al. \cite{AgarwalDutta} have demonstrated the
geotechnical characteristics of recycled aggregates, including C\&D
waste, and reported that these aggregates possess considerable engineering
properties, which can further be improved by stabilization.

Besides, Saride and Jallu \cite{SarideJallu} have shown that about
60\% of RAP can be used in base layers stabilized with low calcium
fly ash geopolymer through a series of laboratory studies and field
test sections. Geopolymer is an inorganic alumina-silicate polymer
material synthesized using alkaline activation of alumina-silicate
source materials \cite{Davidovits}. Several researchers have reported
that the geopolymer stabilized materials possess higher compressive
strength and durability characteristics \cite{Song,Fernandez}. Nevertheless,
if the reactive silica content is low in the source material (fly
ash), in addition to the alkali hydroxides, an external supplement
of silica is suggested to improve the formation of geopolymer structures.
In general, sodium silicate ($Na_{2}SiO_{3}$) or potassium silicate
($K_{2}SiO_{3}$) are supplemented for this purpose.

It was reported that the stabilization of RAP with fly ash geopolymer
would not produce a unique product as it depends on the quality of
fly ash and the RAP \cite{SarideJallu}. They have also demonstrated
that the layer coefficients used in the design of flexible pavements
by the American Association of State Highway and Transportation Officials
(AASHTO) method produce either conservative or unsafe base course
thicknesses when fly ash geopolymers are used. They have proposed
a new set of layer coefficients for FRG stabilized RAP bases.

\subsection{Design of a Flexible Pavement using Reclaimed Material}

Based on the recently developed novel approaches to accommodate secondary/alternative
materials in pavement construction, a full-scale study was undertaken
to evaluate the performance of fly ash stabilized RAP as a base course
for a flexible pavement \cite{SarideChallapalli}. As part of the
study, extensive laboratory and field studies were undertaken.

The key features of the study are outlined below: 
\begin{itemize}
\item As high as eighty percent of virgin aggregate was successfully replaced
with RAP aggregate. About 20\% to 30\% of low calcium fly ash was
utilized by activating the fly ash with alkali activators. About 350
metric tons of fly ash were consumed per lane kilometer of national
highway construction, which amounts to 1400 tons of fly ash required
for 4-lane express highways. 
\item Most of the pavement design guidelines across the world limit the
amount of RAP in the base course up to 30\% by weight of the virgin
aggregates (VA), due to the presence of aged bitumen coating on the
RAP aggregates. It was demonstrated that at least 60\% of the VA could
be replaced while meeting all design requirements. 
\item The critical design parameters, including the Unconfined Compressive
Strength (UCS) and Resilient Modulus (Mr) characteristics of the FA
geopolymer stabilized RAP:VA blends were found to meet the threshold
values proposed by various design codes. 
\item Besides, exposure of these mixes to the severe moisture and temperature
variations may alter the cementation. The permanency of the stabilizer
was also verified through rigorous wet-dry durability and leachate
studies. The UCS tests, durability tests, and resilient modulus test
(Mr) result indicated that the strength loss of RAP: VA mixes were
very minimal and were found suitable for the base course applications. 
\item Based on the laboratory performance of the product, the performance
of the design mix was studied under actual traffic conditions for
three years by laying a trial road stretch (200 m) as part of a state
highway (SH 207) near Vijayawada in Andhra Pradesh. 
\item Overall, based on the experimental and field data, a reduction in
pavement base course thickness by about 30\% and cost of construction
by about 20\% was achieved with a similar or better performance of
the pavement. 
\end{itemize}

\subsection{Variables Influence the Sustainability of Pavements}

\begin{figure}
\includegraphics[scale=0.7]{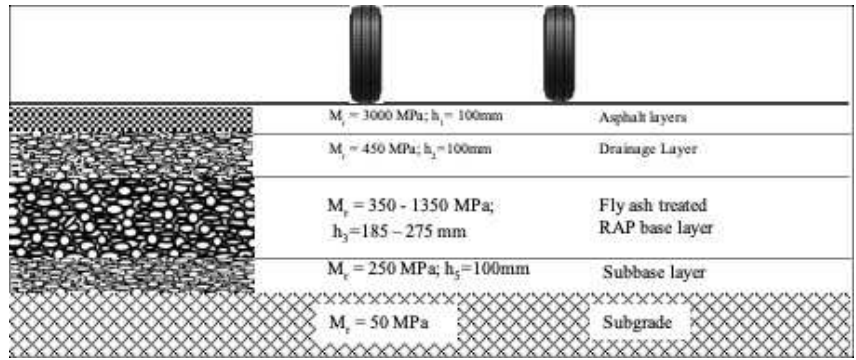}

\caption{\label{Wheels}A typical cross-section of a flexible pavement with
recycled base material}
\end{figure}

\begin{figure}
\begin{description}
\item [{\includegraphics[scale=0.3]{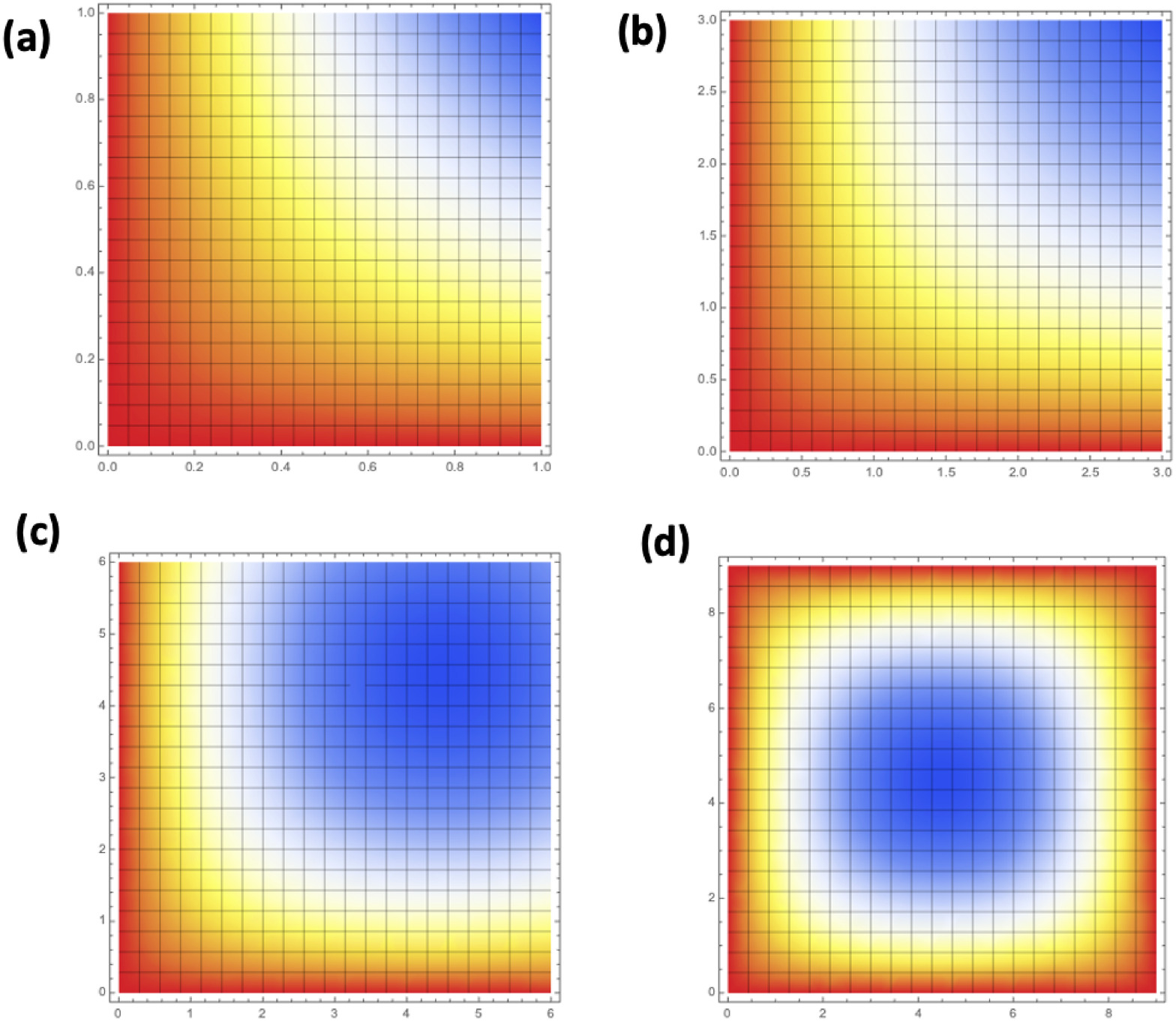}}]~
\end{description}
\caption{\label{fig:Plotting-of-the}Plotting of the model for $k=2$ based
on two variables RAP and VA in the Table \ref{tab:Description}. \textbf{(a)
}$\psi_{1}\rightarrow[0,1]$ and $\psi_{2}\rightarrow[0,1]$, \textbf{(b)
}$\psi_{1}\rightarrow[0,3]$ and $\psi_{2}\rightarrow[0,3]$, \textbf{(c)
}$\psi_{1}\rightarrow[0,6]$ and $\psi_{2}\rightarrow[0,6]$, \textbf{(d)
}$\psi_{1}\rightarrow[0,9]$ and $\psi_{2}\rightarrow[0,9]$. We have
used the initial and boundary conditions as follows: $H(\psi_{1},\psi_{2},0)=0,$
$H(0,\psi_{2},t)=st$, $H(9,\psi_{2},t)=st$, $H(\psi_{1},0,t)=st$,
$H(\psi_{1},9,t)=st$ and $s=10,t\in[0,1000].$}
\end{figure}

\begin{figure}
\begin{description}
\item [{\includegraphics[scale=0.3]{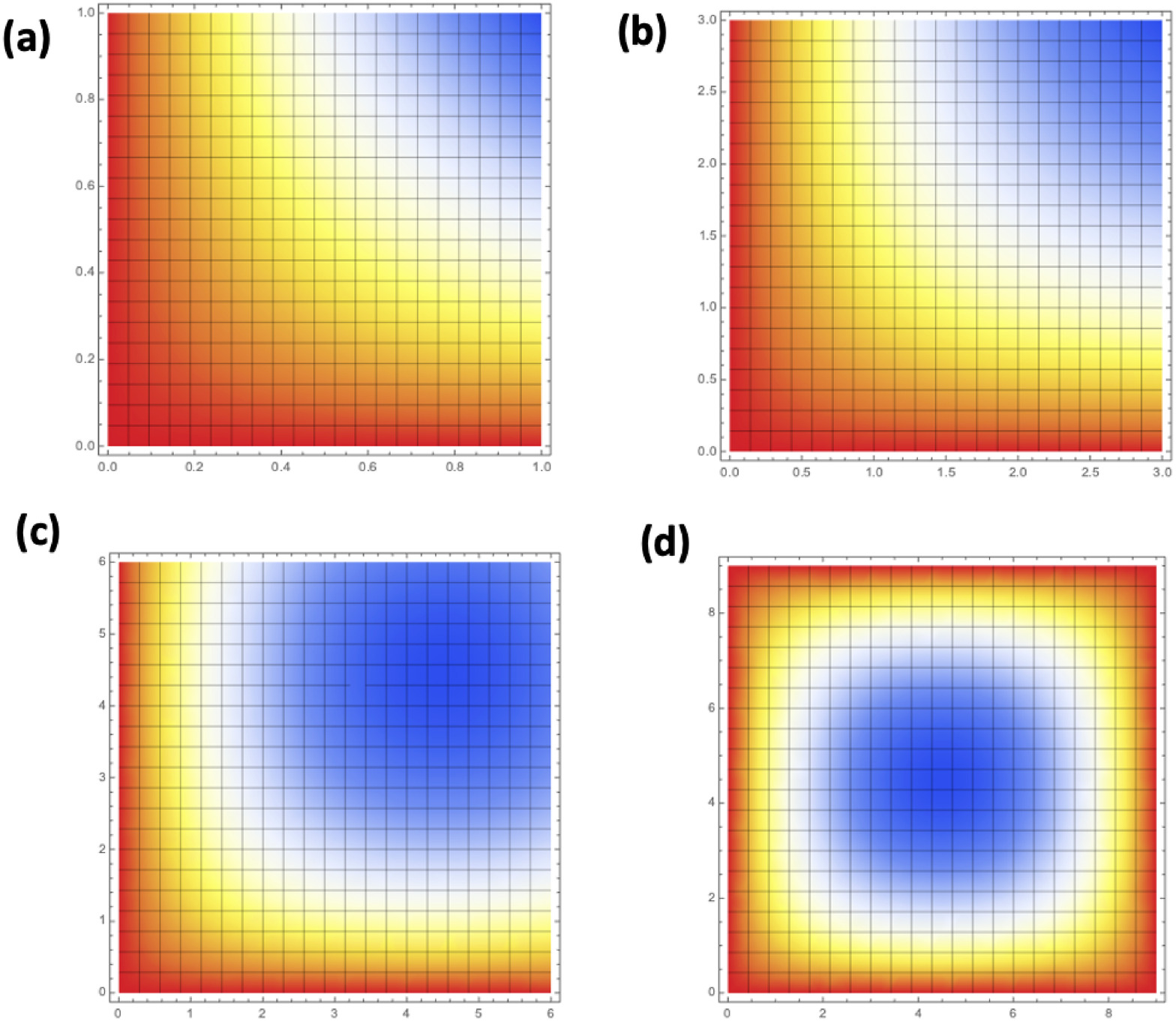}}]~
\end{description}
\caption{\label{fig:Plotting-of-the-density6:4}Plotting of the model for $k=2$
based on two variables RAP and VA in the Table \ref{tab:Description}.
\textbf{(a) }$\psi_{1}\rightarrow[0,4]$ and $\psi_{2}\rightarrow[0,6]$,
\textbf{(b) }$\psi_{1}\rightarrow[0,6]$ and $\psi_{2}\rightarrow[0,9]$,
\textbf{(c) }$\psi_{1}\rightarrow[0,8]$ and $\psi_{2}\rightarrow[0,12]$,
\textbf{(d) }$\psi_{1}\rightarrow[0,10]$ and $\psi_{2}\rightarrow[0,15]$.
We have used the initial and boundary conditions as follows: $H(\psi_{1},\psi_{2},0)=0,$
$H(0,\psi_{2},t)=st$, $H(18,\psi_{2},t)=st$, $H(\psi_{1},0,t)=st$,
$H(\psi_{1},18,t)=st$ and $s=10,t\in[0,1000].$}
\end{figure}

To establish the sustainability index for pavements, a typical cross-section
of flexible pavement is shown in Figure \ref{Wheels} is proposed
to build on a weak sub-grade with a resilient modulus of 30 MPa. The
pavement is designed according to \cite{IndianRoads} for the traffic
of 20 million standard axles (msa). The description of the variables
used and the range are presented in Table \ref{tab:Description}.
According to Indian Roads Congress \cite{IndianRoads}, the standard
base layer thickness, when virgin aggregate (VA) is used, is found
to be 275 mm with an 80 mm thick asphalt layer and 200 mm sub-base
layer (See Table \ref{tab:Description-and-range TableA}). Based on
the above discussion, if the base layer is chosen to be replaced with
a recycled aggregate base (RAP), the thickness of the base layer may
be reduced, as suggested in Table \ref{tab:Description-and-range TableA}.
As discussed above, since VA is replaced with RAP, stabilization is
inevitable to meet the structural strength and stiffness criteria
of the base layer. Further, to enhance sustainability in the construction,
fly ash was adopted against the conventional stabilizer such as cement.
To quantify the effectiveness of the method, the sustainability indicator
for the proposed approach may be evaluated to choose an appropriate
design mix for the pavement design. To calculate sustainability, several
design mixes were considered ranging from replacing 50\% VA with 50\%
RAP stabilized with fly ash. The variation in the combination of RAP:VA
blends are presented in Table \ref{tab:Description-and-range TableA}.
The suitability of each design mix was established by conducting the
strength in terms of unconfined compressive strength (UCS) and resilient
modulus (Mr) \cite{ArulRajah,SarideJallu,Avirneni}. The Mr values
corresponding to each design mix (See Table \ref{tab:Description-and-range TableA})
are used to obtain the design thickness of the base layer for each
combination. It is very clear from the Table \ref{tab:Description-and-range TableA}
that if the VA is replaced by 100\% RAP, the overall thickness of
the pavement has reduced from 555 mm to 510 mm (8\% reduction) if
it is 50\% the reduction is 18\% and 14.5\% for 60\% replacement for
the similar performance of the pavement. However, if the overall performance
in terms of safety and cost is considered, a 60\% replacement of VA
is more sustainable. We have plotted Figure \ref{fig:Plotting-of-the}
and Figure \ref{fig:Plotting-of-the-density6:4} based on the data
collected in this section. The scenario presented in Table \ref{tab:Description}
and Table \ref{tab:Description-and-range TableA} is evaluated using
the proposed sustainability model (\ref{PDE-independent}) by analyzing
the density functions as explained in Figure \ref{fig:Plotting-of-the}
and Figure \ref{fig:Plotting-of-the-density6:4}. Figure \ref{fig:Plotting-of-the}
depicts that when the range is increased from 1 to 9, the relative
influence of the variables is reducing, indicating that when the proportion
of RAP and VA are the same. Figure \ref{fig:Plotting-of-the-density6:4}
is plotted by considering the RAP and VA proportion as 60\% and 40\%,
respectively. Now the influence of RAP on VA is skewed when the range
increases. A similar analysis is suggested to perform by considering
all the variables in the given scenario to estimate the sustainability
index.

\begin{landscape}

\begin{table}
\begin{tabular}{|c|c|c|c|}
\hline 
Sr. No. & Variables & Description & Range of Values\tabularnewline
\hline 
\hline 
1 & RAP & $\begin{array}{c}
\text{Reclaimed asphalt pavement material, }\\
\text{obtained from the milling of distressed pavements}
\end{array}$ & 50, 60, 80\tabularnewline
\hline 
2 & VA & $\text{Natural aggregates used for road construction}$ & 50, 40, 20\tabularnewline
\hline 
3 & FA & $\text{Fly ash, obtained from coal combustion in power plants}$ & 20, 30\tabularnewline
\hline 
4 & Mr & $\begin{array}{c}
\text{Resilient modulus, a stiffness parameter used}\\
\text{in the design of pavement layer thickness}
\end{array}$ & 350 to 1350\tabularnewline
\hline 
\end{tabular}

\caption{\label{tab:Description}Description and range of variables used in
the pavement design}
\end{table}

\begin{table}
\begin{tabular}{|c|c|c|c|c|c|c|c|c|}
\hline 
Sr. No. & Mix Combination & \multirow{1}{*}{AC} & Drainage & Sub-base & Base & $\begin{array}{c}
\text{Total }\\
\text{Thickness}
\end{array}$ & $\begin{array}{c}
\text{Base Mr (MPa) }\\
\text{considered}
\end{array}$ & Reference\tabularnewline
\hline 
\hline 
1 & 0R:100VA & 80 & NA & 200 & 275 & 555 & 350 & IRC:37 (2018)\tabularnewline
\hline 
2 & 50R:50V+20F & 70 & 100 & 100 & 185 & 455 & 1344 & -\tabularnewline
\hline 
3 & 60R:40V+20F & 70 & 100 & 100 & 195 & 465 & 1191 & \cite{SarideAvimeni}\tabularnewline
\hline 
4 & 80R:20V+20F & 70 & 100 & 100 & 205 & 475 & 988 & \cite{Avirneni}\tabularnewline
\hline 
5 & 100R:0VA+20F & 70 & 100 & 100 & 240 & 510 & 565 & \cite{ArulRajah}\tabularnewline
\hline 
6 & 50R:50V+30F & 70 & 100 & 100 & 195 & 465 & 1156 & -\tabularnewline
\hline 
7 & 60R:40V+30F & 70 & 100 & 100 & 205 & 475 & 968 & \cite{SarideJallu}\tabularnewline
\hline 
8 & 80R:20V+30F & 70 & 100 & 100 & 215 & 485 & 824 & \cite{SarideChallapalli}\tabularnewline
\hline 
\end{tabular}

\caption{\label{tab:Description-and-range TableA}Design of a flexible pavement
with recycled materials in the base layer. Note: {*}The pavement is
designed using IRC 37:2018 guidelines for the design traffic of 20
million standard axles (msa) and subgrade modulus of 50 MPa. \textbf{Remarks:
}Subgrade : Mr = 50 MPa, Sub-base Layer: Mr =250 MPa, Drainage Layer:
Mr =450 MPa, AC Layer: Mr =3000 MPa}
\end{table}

\end{landscape}

\section{Conclusions}

\label{sec:5}\textcolor{black}{There is a great need for collaborative
efforts to build mathematical models for measuring and understanding
sustainability based on the perception of the people who work on sustainability
and who work for sustainability. The models (\ref{PDE-independent})
and (\ref{PDE-dependent}) are very flexible, can accommodate a variety
of options to measure sustainability that could arise from the cross-discipliners
team of scientists. Using these models, one can measure the overall
sustainability of human life in a country or a region as a result
of several variables by treating each factor independently (model
(\ref{PDE-independent})) or resultant measure of sustainability as
a result of various overlapping variables (model (\ref{PDE-dependent})).
When there are not enough evidence on the dependencies of variables
for a particular population concerned to study sustainability, we
can start the analysis using the model (\ref{PDE-independent}). In
fact, the model (\ref{PDE-independent}) is a subset of the model
(\ref{PDE-dependent}) under certain situations. Usually, second-order
PDEs arise in most of the natural and engineering sciences situations.
However, higher-order PDEs are also prevalent in handling questions
related to mechanics and elasticity. The degree and size of information
to be used at a country or regional level could be data-dependent,
hence, we have proposed flexible models}\textbf{\textcolor{black}{.
}}\textcolor{black}{Mathematical models might help to understand factors
and features of sustainability, however, for practical solutions,
these models must reach mainstream developmental activities. Obtaining
relevant data from the seven variables mentioned might not be an easy
exercise and could involve time-consuming efforts by people involved
in sustainability research and other key people in political and government
circles. }

\textcolor{black}{The data from atmosphere, oceans, vegetation, food,
wetlands, species, and several environmental parameters could generate
multiple sources of errors and bias. Hence, there is a requirement
for the estimation of variability. As pointed out in \cite{Lev} there
is good scope for statistical thinking as well. There are documented
arguments on how collective efforts by social, political, government
set-up could form a network that could help to maintain adaptability
and transformability in ecological dynamics \cite{Gund} and sustainable
agriculture and energy and impact of climate in developing countries
\cite{Prett,Ravi}. There are several advances in computational techniques,
and facilities to conduct global level high intense computational
experiments of complex mathematical models so that the global community
can easily come together for sustainable development. }

\textcolor{black}{Our approach of modeling the phenomena is also computationally
challenging, in terms of numerical approximation and numerical solutions
to the proposed PDEs. There is no unified approach for computing Riemann-Stieltjes
weights when they are used in PDEs. Even the decision of using an
appropriate partition of the intervals $\left[\omega_{a},\omega_{b}\right]$
needs construction of complex algorithms because we have mixed weights
with the PDEs. In fact, this kind of modeling initiation and developing
global indices will help to strengthen the data collection and information
gathering activities. Governments and non-government agencies like
the UN, the World Bank, MacArthur Foundation, Gates Foundation, etc
needs to encourage cross-disciplinary modeling research teams and
data collection programs. There is an urgent need to bring al the
activities into a common platform and need for standardization of
the sustainability control activities across the globe. Countries
should take efforts to bring transparency in the data collection methods
and definitions of variables for the global health of the human population.
There is a need to start providing annual sustainability indices for
the country, region (formed by a group of countries), continent and,
World. } 

\section*{Acknowledgements}

We thank the referee for constructive suggestions to improve our original
draft. This has helped us to add more data, a numerical example, simulations.
Work for SS is supported by Department of Science and Technology (DST),
New Delhi, India under the grant No. DST/TSG/STS/2013/40.

\end{document}